\documentclass[a4paper,12pt]{amsart}
\usepackage{amsmath, amsthm, amscd, amsfonts}

\newtheorem{theorem}{\quad Theorem}[section]
\theoremstyle{definition}
\newtheorem{definition}[theorem]{\quad Definition}

\theoremstyle{theorem}
\newtheorem{corollary}[theorem]{\quad Corollary}
\newtheorem{remark}[theorem]{\quad Remark}

\newtheorem{proposition}[theorem]{\quad Proposition}
\newtheorem{lemma}[theorem]{\quad Lemma}
\newcommand{\lb}{\label}
\begin{document}

\centerline{}

\centerline{}
\title[Roueentan and Khosravi]{}
\centerline{\Large{\bf On Hopfian(co-Hopfian) and Fitting $S$-acts (I)}}

\centerline{}

\centerline{\Large{\bf }}

\centerline{}

\centerline{\bf {Mohammad Roueentan$^{^{*,1}}$},  Roghaieh Khosravi$^{^2}$}
\centerline{}

\centerline{$^{^1}$ Lamerd Higher Education Center, Shiraz University of Technology}
\centerline{ Lamerd, Iran.}
 \centerline{$^{^2}$ Department of Mathematics, Fasa University}
\centerline{Fasa, Iran.}

\thanks{$^{^*}$Correspondenc: rooeintan@sutech.ac.ir}
\maketitle

\begin{abstract}
 The main purpose of the present work is an investigation of the notions Hopfian (co-Hopfian) acts whose their surjective (injective) endomorphisms are isomorphisms. 
While we investigate conditions that are relevant to these classes of acts, their interrelationship with some other concepts for example quasi-injective and Dedekind-finite acts is studied. Using Hopfian and co-Hopfian concepts, several conditions are given for a quasi-injective act to be Dedekind-finite. Moreover we bring out some properties of strongly Hopfian and strongly co-Hopfian $S$-acts. Ultimately we introduce and study the concept of Fitting acts and over a monoid $S$, some equivalent conditions are found to have all its finitely generated (cyclic) acts Fitting. It is shown that an $S$-act is Fitting if and only if it is both strongly Hopfian and strongly co-Hopfian.  
\end{abstract}

~~~~~~~~~

\section{Introduction and Preliminaries}
The study of Hopfian and co-Hopfian notions for rings and modules was initiated by Hiremath and Varadarajan (\cite{hi}, \cite{va}) and because of their importance has been continued in many articles. The main body of the present work consists in an investigation of the counterpart concepts for acts over monoids.\\  In \cite{fa} Hopfian, co-Hopfian, strongly Hopfian and strongly co-Hopfian acts are defined and some properties of them are investigated. Suppose $A$ is an $S$-act over a monoid $S$, then $A$ is called $Hopfian$ ($co$-$Hopfian$) if any surjective (injective) endomorphism of $A$ is an isomorphism (see \cite{fa}). Also $A$ is $strongly$ $Hopfian$($strongly$ $co$-$Hopfian$) if for any $f\in End(A)$ the ascending (descending) chain $\ker f\subseteq \ker f^2 \subseteq \ker f^3 \subseteq ...$  ($Imf \supseteq Imf^2 \supseteq Imf^3 \supseteq... $) is stationary.
In Section 2, we study Hopfian and co-Hopfian acts presenting basic results and their relationship with some other acts for example quas-injective and Dedekind-finite acts.

In Section 3, while studying the concepts of strongly Hopfian and strongly co-Hopfian $S$-acts, we introduce the concept of Fitting acts and provide useful information cocerning these new classes of acts. Over a monoid $S$ an $S$-act $A$ is called a $Fitting$ act if for any endomorphism $f$ of $A$, there exists a natural number $n$ such that $\ker f^n \vee \mathcal{K}_{{\rm{Im}}f^n}=\bigtriangledown _A$ and $\ker f^n\cap \mathcal{K}_{{\rm{Im}} f^{n}}=\Delta_A$, where $\mathcal{K}_{{\rm{Im}}f^n}=\rho_{{\rm{Im}} f^{n}}$. It is shown that an $S$-act $A$ is Fitting if and only if $A$ is both strongly Hopfian and strongly co-Hopfian. Finally equivalent statements are found on a monoid $S$ to have all its finitely generated (cyclic) acts Fitting.\\

Throughout this article $S$  denotes a monoid and an $S$-act $A_S$  (or $A$) is a unitary $S$-act. Recall that a subact $B$ of an $S$-act $A$ is called $large$ (or essential) in $A$ denoted by $B \subseteq' A$, if any $S$-homomorphism $g: A \longrightarrow C$ such that $g|_B$ is a monomorphism is itself a monomorphism (see \cite{kilp}). From \cite{Rouu}, an $S$-act $A$ is said to be $uniform$, provided that every non-zero subact of $A$ is essential. An element $\theta \in A$ is called a zero element if $\theta s=\theta$ for every $s\in S$. Moreover the one element act is denoted by $\Theta=\lbrace \theta\rbrace$. For an $S$-act $A$, by $E(A)$, we mean the injective envelope $A$. An equivalence relation $\rho$ on an $S$-act $A$ is called a $congruence$ on $A$ if $a~\rho ~b$ implies $(as)~ \rho ~(bs)$ for every $a,b \in A, s \in S$. For an $S$-act $A$ the diagonal relation $\{(a,a)\,|\,a\in A\}$ on $A$ is a congruence on $A$ which is denoted by $\Delta_A$. Also if $B$ is a subact of $A$, then the congruence $(B\times B)\cup \Delta_A$ on $A$ is denoted by $\rho_B$ and is called $Rees$ $congruence$ by the subact $B$. We denote the set of all congruences on $A$ by Con($A$). We encourage the reader to see \cite{kilp} for basic results and definitions relating to acts over monoids.

\section{Hopfian and co-Hopfian $S$-acts}
In this section we investigate preliminary and basic properties of Hopfian and co-Hopfian $S$-acts over a monoid $S$. First we recall some definitions.

\begin{definition}
Let $A$ be a right $S$-act over a monoid $S$. 
\item[{\rm (i)}] (\cite{fa})  $A$ is called $Hopfian$ ($co$-$Hopfian$) if any surjective (injective) endomorphism of $A$ is an automorphism.
\item[{\rm (ii)}]The monoid $S$ is called $right$ $Hopfian$ ($right$ $co$-$Hopfian$) if the $S$-act $S_S$ is Hopfian (co-Hopfian). By a similar definition, if the left $S$-act $_SS$ is Hopfian (co-Hopfian) then $S$ is said to be a $left$ $Hopfian$ ($left$ $co$-$Hopfian$).
\end{definition}

 Recall that for an element $a$ of an $S$-act $A$, the homomomorphism $\lambda_a:S_S \longrightarrow A$ is defined by $\lambda_{a}(s)=as$ for every $s\in S$. Since any endomorphism of $S_S$ is of the form $\lambda_{a}$ for some element $a\in S$, it is easy to see that $S$ is right Hopfian if and only if any right invertible element of $S$ is left cancellative and also $S$ is right co-Hopfian if and only if any left cancellative element of $S$ is right invertible.
  From \cite{kilp}, the trace of an $S$-act $B$ in an $S$-act $A$ is defined by $\displaystyle{tr(B,A):=\bigcup_{\varphi\in Hom(B,A)}  \varphi (B)}$. From \cite{Rou} an $S$-act $A$ is said to be $strongly$ $duo$ if for any subact $B$ of $A, tr(B,A)=B$.

\begin{remark} \label{pr2.1}
\item[{\rm (i)}] By Proposition 1.7 of \cite{Rou}, any stronly duo $S$-act is co-Hopfian.
\item[{\rm (ii)}] By Proposition 3.2 of \cite{Rouu}, any injective unifrom (projective uniform) $S$-act is co-Hopfian (Hopfian).
\item[{\rm (iii)}] It is easy to check that any commutative monoid is right (left) Hopfian.
\end{remark}

From \cite{Rou}, a subact $B$ of an $S$-act $A$ is called $fully$ $invariant$ if for any endomorphism $f$ of $A, f(B)\subseteq B$.

\begin{proposition} \label{pr2.2} Let $S$ be a monoid then the following hold:
\item[{\rm (i)}] Let $B\subseteq A$ be $S$-acts such that $B$ be a fully invariant subact of $A$. If $B$ and $A /B$ are Hopfian (co-Hopfian), then so is $A$.
\item[{\rm (ii)}] Let $\{A_i\}_{i\in I}$ be a family of $S$-acts. If  ${~~\coprod_{i\in I} A_i}$ is Hopfian (co-Hopfian), then each $A_i$ is Hopfian (co-Hopfian).
Also if $~~\prod_{i\in I} A_i$ is Hopfian (co-Hopfian), then each $A_i$  is Hopfian (co-Hopfian).
\item[{\rm (iii)}] If $A$ and $B$ are indecomposable $S$-acts, then $A$ and $B$ are Hopfian (co-Hopfian) if and only if $ A\amalg B$ has the same property.
\end{proposition}

\begin{proof}
(i) First suppose $B$ and $A /B$ are Hopfian and $f: A\longrightarrow A$	 is an epimorphism. Define $g: A/B \longrightarrow A/B$ by $g(\bar{a})= \overline{f(a)}$ for any element $a\in A$. It is clear that $g$ is an epimorphism which is an isomorphism by assumption. It implies that the homomorphism $h:= f|_B :B\longrightarrow B $ is an epimorphism and so by assumption is an  isomorphism. Now we can easily see that $f$ is a monomorphism. By a similar proof, we can prove the co-Hopficity case. \\
(ii) Trivial.\\
(iii) We will only regard the co-Hopficity case and the other case is similar. The sufficient part holds by part (ii). Now suppose $A$ and $B$ are co-Hopfian and $ f:A\amalg B \longrightarrow A\amalg B$ is a monomorphism. Since $A$ and $B$ are indecomposabe, three cases can be considered.\\
(a) $\bold {f(A)\subseteq A}$ and $\bold {f(B)\subseteq B}$: The proof is clear.\\
(b) $\bold {f(A)\subseteq B}$ and $\bold {f(B)\subseteq A}$: The homomorphisms $fof:A\longrightarrow A$ and $fof:B\longrightarrow B$ are monomorphisms, which are isomorphisms by assumption and the proof is evident.\\
(c) $\bold {f(A), f(B)\subseteq B}$ or $\bold {f(A), f(B)\subseteq A}$: In each of these cases, using the restrictions of $f$ to $A$ and $B$, we have $A=B$ which is a contradiction.
\end{proof}

 The following lemma is a direct consequence of Proposition 3.4 of \cite{hin}. From \cite{Sa}, an $S$-act $A$ is said to be $torsion$ $free$ if for any $s\in S$ and any elements $a,b\in A$, the equality $as=bs$ implies that $a=b$.
\begin{lemma} \label{le2.7}
 Suppose $A$ is a weakly injective $S$-act with a zero and $E(A)$ is a torsion free $S$-act. Then $A$ is injective.
\end{lemma}

Recall that a monoid $S$ is called \textit{right hereditary} if all  right ideals of $S$ are projective.
\begin{proposition} \label{pr2.7}
Let $A$ be an $S$-act with a zero over a right hereditary monoid $S$. If $E(A)$ is a torsion free Hopfian $S$-act, then so is $A$.
\end{proposition}

\begin{proof}
Suppose $f:A\longrightarrow A$ is an epimorphism. Thus there exists a homomorphism $g: E(A)\longrightarrow E(A)$ which is an extension of $f$. Since $S$ is hereditary, by \cite[Theorem 4.11.20]{kilp}, ${\rm{Im}} g$ is a weakly injective $S$-act and by the above lemma is injective. Since $A=f(A)\subseteq {\rm{Im}} g\subseteq E(A)$ and $A\subseteq' E(A)$, we have ${\rm{Im}} g\subseteq' E(A)$. Hence  ${\rm{Im}} g=E(A)$ and so $g$ is an epimorphism. Now, by assumption we infer that $g$ is a monomorphism which implies that $f$ is also a monomorphism as desired.
\end{proof}

\begin{proposition} \lb{pr2.77}
Suppose $A$ is an $S$-act such that for any its non-zero element $a, A\hookrightarrow aS$. Then $A$ is uniform if and only if $E(A)$ is co-Hopfian.
\end{proposition}
\begin{proof}
If $E(A)$ is co-Hopfian and $B$ is a nonzero subact of $A$, then there exists a monomorphim $f: E(A)\longrightarrow E(B)$ which is an isomorphism by assumption. Thus $E(B)$ is also co-Hopfian and $f|_{E(B)}$ is also an isomorphism. Now if $a\in E(A)$, then $f(a)=f(b)$ for some $b\in E(B)$ and since $f$ is a monomorphism, $a=b\in E(B)$. Thus $E(A)=E(B)$ and by \cite{Rouu}, $B$ is a large subact of $A$ which implies that $A$ is uniform. The converse holds by Lemma {2.11} and Proposition {3.2} of \cite{Rouu}.
\end{proof}

\begin{definition}
Let $A$ be an $S$-act over a monoid $S$. Then $A$ is called $quasi$-$retractable$ ($quasi$-$coretractable$) if any surjective (injective) endomorphism of $A$ is right invertible (left invertible).
\end{definition}

From \cite{Ahs} recall that an $S$-act $A$ is $quasi$-$projective$ if for every $S$-act $C$, every $S$-homomorphism $f:A\longrightarrow C$ can be lifted with respect to every epimorphism $g:A\longrightarrow C$, i.e., there exists a homomorphism $h:A\longrightarrow A$ such that $f=gh$. Also From \cite{Sa} an $S$-act $A$ is called \textit{quasi-injective} if it is injective relative to all inclusions from its subacts, i.e., for any homomorphism $f: B \longrightarrow A$, there exists a homomorphism $h:A\longrightarrow A$ such that $h|_B=f$ where $B$ is a subact of $A$. \\
We can check that any quasi-projrctive (quasi-injective) $S$-act is quasi-retractable (quasi-coretractable). Also any Hopfian (co-Hopfian) $S$-act is quasi-retractable (quasi-coretractable).\\
The proof of the following proposition is straightforward and will be omited.

\begin{proposition} \label{pr2.3} Let $S$ be a monoid then the following hold:
\item[{\rm (i)}] Any quasi-retractable co-Hopfian $S$-act is Hopfian.
\item[{\rm (ii)}] Any quasi-coretractable Hopfian $S$-act is co-Hopfian.
\item[{\rm (iii)}] If $A$ is a Hopfian (co-Hopfian) $S$-act, then every left or right invertible endomorphism of $A$ is invertible.
\end{proposition}

By the above proposition, we conclude that any quasi-projective co-Hopfian $S$-act is Hopfian and any quasi-injective Hopfian $S$-act is co-Hopfian.\\

\begin{definition}
Let $A$ be an $S$-act. Then $A$ is called $mono$-$uniform$ if for any monomorphism $f:A\longrightarrow A, f(A)\subseteq' A$.
\end{definition}
 It is clear that any uniform (co-Hopfian) $S$-act is mono-unifrom. Also the right $S$-act $S_S$ is mono-uniform if and only if for any left cancellative element $x\in S, xS\subseteq' S$.
\begin{proposition} \label{pr2.33} Let $S$ be a monoid then the following hold:
\item[{\rm (i)}]  If the right $S$-act $S_S$ is mono-uniform, then $S$ is right Hopfian.
\item[{\rm (ii)}]  If the right $S$-act $S_S$ is mono-uniform and $E(S_S)$ is a cyclic $S$-act, then $S_S$ is injective.
\end{proposition}
\begin{proof}
(i) If $f:S_S\longrightarrow S_S$ is an epimorphism, then for some elements $x,s\in S, f(1)=s$ and $1=f(x)$. Thus $1=f(x)=sx$ and so the homomorphism $\lambda_{x}:S_S \longrightarrow S_S$ is a monomorphism. Thus by assumption $xS=\lambda_{x}(S_S)\subseteq'S$ which implies that $f$ is a monomorphism.\\
(ii) Since $E(S_S)$ is a cyclic $S$-act there exists an epimorphism $f:S_S \longrightarrow E(S_S)$. Now, by projectivity of $S_S$, there exists a homomorphism $g:S_S \longrightarrow S_S$ such that $fg=1_S$. We imply that $g$ is a monomorphism and $\mathcal{K}_{{\rm{Im}}g}\cap \ker f=\Delta_S$. Since $\mathcal{K}_{{\rm{Im}}g}\subseteq' S, \ker f=\Delta_S $ and hence $f$ is a monomrphism. Thus $f$ is an isomorphism and the proof is completed.
\end{proof}
 Since any quasi-injective Hopfian $S$-act is co-Hopfian, by the above proposition we have the following corollary.
\begin{corollary} \label{co2.34}
Suppose $S$ is a monoid and $E(S_S)$ is a cyclic $S$-act. Then $S_S$ is mono-uniform if and only if $S$ is right co-Hopfian.
\end{corollary}

From \cite{kh}, an $S$-act $A$ is called $noetherian$ ($artinian$) if $Con(A)$ satisfies the ascending (descending) chain condition. Also an $S$-act $A$ is called $Rees$ $noetherian$ ($Rees$ $artinian$) if it satisfies the ascending (descending) chain condition on its subacts. For the sake of simplicity, we denote " ascending chain condition" and " descending chain condition" by "a.c.c" and "d.c.c ", respectively.

\begin{lemma} \label{le2.11}
Suppose $A$ is a Rees noetherian $S$-act with a zero and $B$ is a subact of $A$ such that $A\cong A/B$, then $B=\Theta$.	
\end{lemma}
\begin{proof}
Since $A$ is Rees noetherian, by Proposition 2.1 of \cite{kh}, we can suppose $B$ is a maximal subact of $A$ with respect to the condition $A\cong A/B$. If $\varphi:A \longrightarrow  A/ B$ is an isomorphism, then obviously
 $ A / B \cong  {(A/B)} /{\varphi (B)}$	and so $A\cong {(A/B)} /{\varphi (B)}$. If $ B$ is non-zero then $\varphi (B) $ is a non-zero subact of $A/B$. Thus there exists a subact $C$ of $A$ which contains $B$ properly and $A\cong A/C$, a contradiction to the maximality of $B$. Therefor $B=\Theta$.
\end{proof}

\begin{proposition} \label{co2.12}
Suppose $A$ is a quasi-injective co-Hopfian $S$-act with a zero. If $B$ is an essential subact of $A$ such that $A/B$ is Rees noetherian, then $B$ is co-Hopfian.
\end{proposition}
\begin{proof}
By way of contradiction suppose that $f:B \longrightarrow B $ is a monomorphism and $f(B)\neq B$. Therefor since $B\subseteq' A$, using quasi-injectivity of $A$, $f$ can be extended to an automorphism $\bar{f}$ of $A$. Thus $ A / B \cong  {\bar{f}(A)} /{\bar{f} (B)}= A / f(B)$. Hence $A / B \cong   A / f(B)$ and by assumption $A / f(B)$ is Rees noetherian. Thus $\mathcal{B}= B/f(B)$ is a non-zero subact of $ \mathcal{A}= A / f(B)$ such that $\mathcal{A}/\mathcal{B}$ is Rees noetherian and by the above lemma, we have a contradiction.
\end{proof}

\begin{corollary}\label{co2.13}
Suppose for an $S$-act $A$, $E(A)$ is co-Hopfian and $E(A)/A$ is Rees noetherian, then $A$ is co-Hopfian.
\end{corollary}

\begin{proposition} \label{pr2.14}
Suppose $\mathcal P$ is a property of $S$-acts which is preserved under isomorphism and $A$ is an $S$-act which has this property. If $A$ satisfies a.c.c on congruences $\lambda\in Con(A)$ such that $A/ \lambda$ has the property $\mathcal P$, then $A$ is Hopfian.
\end{proposition}

\begin{proof}
By way of contradiction, suppose $A$ is not Hopfian. Thus there exists a non-diagonal congruences $\sigma_1\in Con(A)$ such that $A\cong A/ \sigma_1$. Thus $A/ \sigma_1$ is an $S$-act which has the property and is not Hopfian. By a similar way we can find a non-diagonal congruences $\sigma_2\in Con(A)$ such that $\sigma_1\subset \sigma_2$. Repeating this process yields an ascending chain of congruences of the form $\sigma_1\subset \sigma_2\subset \sigma_3\subset... $ such that for any $i=1,2,3,.. ~,A/ \sigma_i $ has the property $\mathcal P$, a contradiction.
\end{proof}

By the previous proposition if we consider co-Hopficity as the property $\mathcal P$, then we have the following corollary.
\begin{corollary} \label{pr2.15}
 Suppose $A$ is an $S$-act which satisfies a.c.c on congruences $\lambda\in Con(A)$ such that $A/ \lambda$ is co-Hopfian, then $A$ is Hopfian.
\end{corollary}

\begin{corollary} \label{pr2.16}
 Suppose $A$ is an $S$-act which satisfies a.c.c on congruences $\lambda\in Con(A)$ such that $A/ \lambda$ is not Hopfian, then $A$ is Hopfian.
\end{corollary}
\begin{proof}
By Propositin \ref{pr2.14} if $\mathcal P$ is the property of being not Hopfian, then the result follows.
\end{proof}

\begin{proposition} \label{pr2.17}
Suppose $\mathcal P$ is a property of $S$-acts which is preserved under isomorphism and $A$ is an $S$-act which has this property. If $A$ satisfies d.c.c on subacts with property $\mathcal P$, then $A$ is co-Hopfian.
\end{proposition}
\begin{proof}
If $A$ is not co-Hopfian, then by a routin argument, we can obtain a strictly descending chain of proper subacts with property $\mathcal P$ that cntradicts our hypothesis.
\end{proof}

\begin{corollary} \label{pr2.18}
Suppose $A$ is Hopfian $S$-act and saisfies d.c.c on Hopfian subacts, then $A$ is co-Hopfian.
\end{corollary}

\begin{corollary} \label{pr2.19}
 Let $A$ be an $S$-act which has d.c.c on its non-co-Hopfian subacts, then $A$ is co-Hopfian.
\end{corollary}
\begin{proof}
By Propositin \ref{pr2.17} it is sufficient to let $\mathcal P$ the property of being non-co-Hopfian.
\end{proof}

\begin{proposition} \label{pr2.20}
Suppose $A$ is an $S$-act which satisfies d.c.c on its nonessential subacts, then $A$ is mono-unifrom.
\end{proposition}
\begin{proof}
The proof is similar to the proof of Proposition \ref{pr2.17}.
\end{proof}

 Concluding this section, we concentrate on the category $Act_0$--$S$ and from now to the end of this section, we suppose that any $S$-act comtained in $Act_0$--$S$. From \cite{Ah}, in the category $Act_0$--$S$, an $S$-act $A$ is called $Dedekind$-$finite$ if $B = \Theta$ is the only S-act for which $A\coprod B\cong A$. Note that in the category $Act_0$--$S$,  $A_1\coprod A_2=A_1\cup A_2$, where $ A_1\cap A_2= \Theta$, and each $A_i$ is called a $0$-direct summand of $A$.  By part (iii) of the previous proposition and Theorem 2.15 of \cite{Ah}, we have the following corollary.
\begin{corollary} \label{co2.4}
Any Hopfian (co-Hopfian) $S$-act is Dedekind-finite.
\end{corollary}

Since any quasi-injective $S$-act is a quasi-coretractable $S$-act, the following theorem is a good improvement of Theorem 2.7 of \cite{Ah}.
\begin{theorem}\label{th2.5}
Suppopse $A$ is a torsion free quasi-coretractable $S$-act. Then the following conditions are equivalent:
\item[{\rm (i)}] $A$ is Dedekind-finite.
\item[{\rm (ii)}] $A$ is co-Hopfian.
\item[{\rm (iii)}] $A$ is mono-uniform.
 \item[{\rm (iv)}] Injective endomorphisms of $A$ map essential subacts to essential subacts.

\end{theorem}
\begin{proof}
(i)$\Rightarrow$ (ii) Suppose $f:A\longrightarrow A$ is a monomorphism. Hence by assumption there exists an endomorphism $g$ of $A$ such that $gf=1_A$. Thus by Theorem 2.15 of
\cite{Ah}, $fg=1_A$ and so $f$ is an epimorphism.\\
(ii)$\Rightarrow$ (iii) Trivial.\\
(iii)$\Rightarrow$ (i) If $g:A\amalg B \longrightarrow A $ is an isomorphism for some $S$-act $B$, then by assumption $ g(A\amalg B)=A\subseteq' A\amalg B$ which implies that $B=\Theta$.\\
(iii)$\Rightarrow$ (iv) If $B$ is an essential subact of $A$ and $f$ is an injective endomorphism of $A$, then by a routine argument, we can see that $f(B) \subseteq' f(A)$ and by assumption the result follows.\\
(iv)$\Rightarrow$ (iii) Clear.\\
\end{proof}
 The next  crollary is obtained from the above theorem.
\begin{corollary} \label{co2.14}
Let $A$ be a torsion free $S$-act. Then $A$ is co-Hopfian if and only if $A$ is mono-uniform and quasi-coretractable.
\end{corollary}

\begin{corollary}\label{co2.55}
Suppopse $A$ is a torsion free quasi-injective $S$-act. Then the following conditions are equivalent:
\item[{\rm (i)}] $A$ is Dedekind-finite.
\item[{\rm (ii)}] $A$ is co-Hopfian.
\item[{\rm (iii)}] $A$ is mono-uniform.
 \item[{\rm (iv)}] Injective endomorphisms of $A$ map essential subacts to essential subacts.
 \item[{\rm (v)}] Any essential fully invariant subact of $A$ is co-Hopfian.
\item[{\rm (vi)}] Any essential fully invariant subact of $A$ is Dedekind-finite .

\end{corollary}

\begin{proof}
The implications (i)$\Leftrightarrow$ (ii)$\Leftrightarrow$ (iii)$\Leftrightarrow$ (iv) hold by Theorem \ref{th2.5}. \\
(v)$\Rightarrow$ (ii) Clear.\\
(ii)$\Rightarrow$ (v) We show that any essential fully invariant subact of $A$ is mono-uniform. First note any fully invariant subact of any quasi injective $S$-act is quasi-injective. Suppose $B$ is an essential fully invariant subact of $A$ and $C\subseteq' B$. Since $B\subseteq' A$ and $A$ is quasi-injective, any injective endomorphism $f$ of $B$ can be extended to an injective endomorphism $g$ of $A$. Now due to equivalence (ii)$\Leftrightarrow (iii),  g(C)\subseteq' A$ and so $g(C)=f(C)\subseteq' B\subseteq A$, as desired.\\
(i)$\Leftrightarrow$ (vi) Duo to the equivalence of the previous conditions can be proved.
\end{proof}

It can be easily checked that for an $S$-act $A$, if a monomorphism $f:A\longrightarrow A$ can be extended to an isomorphism $\bar{f}:E(A)\longrightarrow E(A)$, then $f(A)\subseteq'A$.

\begin{corollary}\label{co2.9}
Suppose $A$ is a quasi-injective $S$-act with the same conditions in the previous theorem. If $E(A)$ is co-Hopfian then $A$ is co-Hopfian.
\end{corollary}
\begin{proof}
By assumption any monomorphism $f:A \longrightarrow A$ can be extended to an isomorphism $\bar{f}: E(A)\longrightarrow E(A)$. Since $A\subseteq' E(A)$ and $\bar{f}$ is an isomorphism, $\bar{f}(A)\subseteq' E(A)$. Moreover $\bar{f}(A)=f(A)\subseteq A$. Thus $f(A)\subseteq'A$ and by Theorem \ref{th2.5} the proof is completed.	
\end{proof}

By a similar proof of the above theorem, we have the next theorem.
\begin{theorem}\label{te2.6}
Suppopse $A$ is a torsion free quasi-retractable $S$-act. Then the following conditions are equivalent:
\item[{\rm (i)}] $A$ is Dedekind-finite.
\item[{\rm (ii)}] $A$ is Hopfian.
\end{theorem}

\begin{corollary}\label{co2.66}
Suppopse $A$ is a torsion free quasi-projective $S$-act. Then the following conditions are equivalent:
\item[{\rm (i)}] $A$ is Dedekind-finite.
\item[{\rm (ii)}] $A$ is Hopfian.
\end{corollary}

\section{On Fitting $S$-acts }
This section deals with the concept of Fitting $S$-acts. First some properties of strongly Hopfian and strongly co-Hopfian acts are studied and then the notion of Fitting $S$-acts is introduced and studied.\\
\begin{definition}
 (\cite{fa}) An $S$-act $A$ is said to be $strongly$ $Hopfian$ ($strongly$ $co$-$Hopfian$) if for any $f\in End(A)$ the ascending (descending) chain $\ker  f\subseteq \ker  f^2 \subseteq \ker  f^3 \subseteq ...$  (${\rm{Im}}f \supseteq {\rm{Im}}f^2 \supseteq {\rm{Im}}f^3 \supseteq... $) is stationary.
\end{definition}
 If $S$ is a monoid, then obviously the right $S$-act $S_S$ is strongly Hopfian if and only if for any element $x\in S$, there exists a natural number $n$ such that for any $s,t\in S$ the equality $x^{n+1}s=x^{n+1}t$ implies that $x^{n}s=x^{n}t$. Also the right $S$-act $S_S$ is strongly co-Hopfian if and only if for any elements $z,x\in S$, there exists an element $t\in S$ and also there exists a natural number $n$ such that $x^{n}z=x^{n+1}t$.

The following proposition is recalled from \cite{fa} .

\begin{proposition} \label{pr3.1}
 Let $A$ be an $S$-act over a monoid $S$. The following conditions are equivalent:
\item[{\rm (i)}] $A$ is strongly Hopfian.
\item[{\rm (ii)}] For any $f\in End(A)$ there exists a natural number $n$ such that $\ker  f^n= \ker  f^{n+1}$.
\item[{\rm (iii)}] For any $f\in End(A)$ there exists a natural number $n$ such that $\ker  f^n\cap \mathcal{K}_{{\rm{Im}} f^{n}}=\Delta_A$.
\end{proposition}

 The following proposition is an improvement of Proposition 3.3 of \cite{fa}. For an $S$-act $A$ by $\bigtriangledown _A$, we mean the congruence $A\times A$.

\begin{proposition} \label{pr3.2}
 Let $A$ be an $S$-act over a monoid $S$. The following conditions are equivalent:
\item[{\rm (i)}] $A$ is strongly co-Hopfian.
\item[{\rm (ii)}] For any $f\in End(A)$ there exists a natural number $n$ such that ${\rm{Im}} f^n={\rm{Im}} f^{n+1}$.
\item[{\rm (iii)}] For any $f\in End(A)$ there exists a natural number $n$ such that $\ker  f^n \vee \mathcal{K}_{{\rm{Im}} f^{n}}=\bigtriangledown _A$.
\end{proposition}
\begin{proof}
(i)$\Leftrightarrow$ (ii) Clear.\\
(ii) $\Rightarrow$ (iii) By a similar proof of Proposition 2.4 of \cite{kh}, can be proved.
(iii) $\Rightarrow$ (ii) Suppose $f\in End(A)$ and for a natural number $n, \ker  f^n \vee \mathcal{K}_{{\rm{Im}} f^{n}}=\bigtriangledown _A$. Clearly ${\rm{Im}} f^{n+1}\subseteq {\rm{Im}} f^{n}$. Thus suppose $y=f^n(x)\in {\rm{Im}} f^n$. Since $(x,y)\in \bigtriangledown _A=\ker  f^n \vee\mathcal{K}_{{\rm{Im}} f^{n}}$, for some natural number $m$, there exist elements $b_1,b_2,...,b_m\in A$ such that $(x,b_1)\in\ \lambda_1, ~(b_1,b_2)\in\ \lambda_2 ....,~ (b_m,y)\in\ \lambda_m$ where for any $i\in\lbrace {1,2,...,m}\rbrace, \lambda_i\in \lbrace {\ker  f^n, \mathcal{K}_{{\rm{Im}} f^{n}}}\rbrace$. If $(x,b_1)\in\mathcal{K}_{{\rm{Im}} f^{n}}$, then $y\in {\rm{Im}} f^{2n}\subseteq {\rm{Im}} f^{n+1}$. Thus suppose $(x,b_1)\in \ker  f^n$ and so $f^n(x)=f^n(b_1)$. Again if $(b_1,b_2)\in\mathcal{K}_{{\rm{Im}} f^{n}}$, then the equality $f^n(x)=f^n(b_1)$ implies the result and otherwise we have $(x,b_2)\in \ker  f^n$. Thus by repeating this method, we get the desired result.
\end{proof}

Considering the above propositions, we can reach the following implications:\\

      Noetherian $\Rightarrow $ strongly Hopfian $\Rightarrow $ Hopfian.

      Rees Artinian $\Rightarrow $ strongly co-Hopfian $\Rightarrow $ co-Hopfian.

\begin{definition}
Let $A$ be an $S$-act. Then we say that $A$ is a $Fitting$ act if for any endomorphism $f$ of $A$ there exists a natural number $n$ such that $\ker  f^n \vee \mathcal{K}_{{\rm{Im}} f^{n}}=\bigtriangledown _A$ and $\ker  f^n\cap \mathcal{K}_{{\rm{Im}} f^{n}}=\Delta_A$.
\end{definition}

 Immediately by Propositions \ref{pr3.1} and \ref{pr3.2}, we have the following theorem.

\begin{theorem}\label{te3.3}
Let $A$ be an $S$-act over a monoid $S$. Then $A$ is a Fitting $S$-act if and only if $A$ is both strongly Hopfian and strongly co-Hopfian.
\end{theorem}

\begin{proposition} \label{pr3.33}
Any noetherian Rees artinain $S$-act is Fittng.	
\end{proposition}
\begin{proof}
Holds by Proposition 2.5 of \cite{kh}.
\end{proof}

\begin{proposition} \label{pr3.4}
Let $A$ be an $S$-act. Then the following hold:
\item[{\rm (i)}] If $A$ is quasi-projective, then $A$ is strongly co-Hopfian if and only if for any $f\in End(A)$ there exists a natural number $n$ and also there exists an endomorphism $\gamma$ of $A$ such that $f^n=f^{n+1}\gamma$.

\item[{\rm (ii)}]  If $A$ is quasi-injective, then $A$ is strongly Hopfian if and only if for any $f\in End(A)$ there exists a natural number $n$ and also there exists an endomorphism $\gamma$ of $A$ such that $f^n=\gamma f^{n+1}$.
\end{proposition}

\begin{proof}
(i) Necessity. Let $f\in End(A)$. By assumption and Proposition \ref{pr3.2}, there exists a natural number $n$ such that ${\rm{Im}} f^n={\rm{Im}} f^{n+1}$. Define $g:A \longrightarrow {\rm{Im}} f^n={\rm{Im}} f^{n+1}$ by $g(x)=f^{n+1}(x)$ and $h:A \longrightarrow {{\rm{Im}}} f^n={\rm{Im}} f^{n+1}$ by $h(x)=f^{n}(x)$ for any $x\in A$. Thanks to the quasi-projectity of $A$, there exists an endomorphism $\gamma$ of $A$ such that $h=g\gamma$.\\
Sufficiency. It is clear by Proposition \ref{pr3.2}.\\
(ii) The proof is similar to that of part (i) by using Proposition \ref{pr3.1}.
\end{proof}
If $A$ is a cyclic $S$-act over a commutative monoid $S$, then $End(A)$ is commutative. Now, by Theorem \ref {te3.3} and the previous proposition, we have the following corollary.

\begin{corollary} \label{pr3.44}
Let $S$ be a commutative monoid. Then the following hold:
\item[{\rm (i)}] Any cyclic quasi-injective strongly Hopfian (quasi-projective strongly co-Hopfian) $S$-act is strongly co-Hopfian (strongly Hopfian)
\item[{\rm (ii)}] Any cyclic quasi-injective strongly Hopfian (cyclic quasi-projective strongly co-Hopfian) $S$-act is a Fitting $S$-act.
\end{corollary}

\begin{corollary} \label{co3.444}
Let $S$ be a monoid. Then the following conditions are equivalent:
\item[{\rm (i)}] Any (finitely generated, cyclic) $S$-act is strongly co-Hopfian..
\item[{\rm (ii)}] Any (finitely generated, cyclic) quasi-retractable $S$-act is strongly co-Hopfian.
\item[{\rm (iii)}] Any (finitely generated, cyclic) quasi-projrctive $S$-act is strongly co-Hopfian.
\item[{\rm (iv)}] Any (finitely generated, cyclic) projective $S$-act is strongly co-Hopfian.
\end{corollary}
\begin{proof}
(iv) $\Rightarrow$ (i) Suppose $A$ is an $S$-act and $f\in End(A)$. If $P$ is a projective $S$-act and $g:P\longrightarrow A$ is an epimorphism, then by projectivity of $P$, there exists $h\in End (P)$ such that $gh=fg$. By assumption and Proposition \ref{pr3.4} there exists a natural number $n$ such that $h^{n}=h^{n+1}{\gamma}$ where $\gamma\in End(P)$. Thus $f^n(A)=f^n(g(P))=gh^n(P)=gh^{n+1}\gamma(P)\subseteq gh^{n+1}(P)=f^{n+1}(g(P))=f^{n+1}(A)$. Thus by Proposition \ref{pr3.2} the result is evident.\\
The other implications are clear. Clearly we can consider the proof for finitely generated and cyclic $S$-acts.
\end{proof}

\begin{proposition} \label{pr3.5}
 Let $A$ be a finitely generated $S$-act over a monoid $S$. Then the following conditions are equivalent:
\item[{\rm (i)}] All factor acts of $A$ are co-Hopfian.
\item[{\rm (ii)}] All factor acts of $A$ are strongly co-Hopfian.
\end{proposition}
\begin{proof}
(i)$\Rightarrow $(ii) For any $\rho\in Con(A)$, every factor of the act $A/ \rho$ is again a factor of $A$ and so it is sufficient to show that $A$ is strongly co-Hopfian when any factor of $A$ is co-Hopfian. Suppose $f:A\longrightarrow A$ is an endomorphism of $A$ and let $\sigma{=\bigcup_{n\geq{1}} \lambda_n}$ where for any ${n\geq{1}}, {\lambda_n}=\ker f^n$. It is clear that $\sigma\in Con(A)$ and for any $(a,b)\in\sigma, (f(a),f(b))\in\sigma$. Define $g:{A/\sigma} \longrightarrow {A/\sigma}$ by $g(\bar{a})=\overline{f(a)}$ for any $a\in A$.
It is easy to see that $g$ is well-defined and also is a monomorphism. Thus by assumption $g$ is an epimorphism. Hence if for some natural $n, \lbrace a_1, a_2,...,a_n\rbrace$ is a generating set for $A$, then there exist elements $\bar{b_1}, \bar{b_2},...,\bar{b_n}$ of $A/\sigma$ such that $\bar{a_i}=g(\bar{b_i})=\overline {f(b_i)}$, for $i=1,2,...,n$ and consequently $(a_i,f(b_i))\in\sigma$ for $i=1,2,...,n$. Therefore by definition of $\sigma$, we can find natural numbers $k_1,k_2,...,k_n$ such that $f^{k_i}(a_i)=f^{{k_i}+{1}}(b_i)$, for $i=1,2,...,n$. Thus if $k=max{\lbrace k_1,k_2,...,k_n\rbrace}$, then ${\rm{Im}}f^k={\rm{Im}}f^{k+1}$ and the proof is completed.\\
(ii)$\Rightarrow $(i) Trivial.
\end{proof}

\begin{corollary} \label{co3.6}
Let $S$ be a monoid. Then the following conditions are equivalent:
\item[{\rm (i)}] All finitely generated $S$-acts are co-Hopfian.
\item[{\rm (ii)}]All finitely generated $S$-acts are strongly co-Hopfian.
\item[{\rm (iii)}] All finitely generated quasi-retractable $S$-acts are strongly co-Hopfian.
\item[{\rm (iv)}] All finitely generated quasi-projective $S$-acts are strongly co-Hopfian.
\item[{\rm (v)}] All finitely generated projective $S$-acts are strongly co-Hopfian.

\end{corollary}
\begin{proof}
(i)$\Leftrightarrow $(ii) Clear by the previous proposition.\\
The other implications hold by Corollary \ref{co3.444}.

\end{proof}

\begin{corollary} \label{co3.7}
Let $S$ be a monoid. Then the following conditions are equivalent:
\item[{\rm (i)}] All cyclic $S$-acts are co-Hopfian.
\item[{\rm (ii)}] All cyclic $S$-acts are strongly co-Hopfian.
\item[{\rm (iii)}] All cyclic quasi-retractable $S$-acts are strongly co-Hopfian.
\item[{\rm (iv)}] All cyclic quasi-projective $S$-acts are strongly co-Hopfian.
\item[{\rm (v)}] All cyclic projective $S$-acts are strongly co-Hopfian.
\item[{\rm (vi)}] $S_S$ is strongly co-Hopfian.
\item[{\rm (vii)}] For any elements $z,x\in S$, there exists an element $t\in S$ and also there exists a natural number $n$ such that $x^{n}z=x^{n+1}t$.
\end{corollary}
\begin{proof}
(i)$\Leftrightarrow $(ii) Clear by the previous proposition. \\
(vi)$\Leftrightarrow $(vii) Notice that any endomorphism $f$ of $S_S$ is equal to a homomorphism $\lambda_x$ for some $x\in S$. \\
(ii)$\Rightarrow$ (vi) Clear.\\
(vi)$\Rightarrow$ (ii) Let $A=aS$ is a cyclic $S$-act and $f$ is an endomorphism of $A$. If $f(a)=ax$ for some $x\in S$ then by assumption and Proposition \ref{pr3.2}, for the homomorphism $g=\lambda_x:S_S\longrightarrow S_S$ there exists a natural number $n$ such that $ {\rm{Im}}g^n={\rm{Im}}g^{n+1}$. Thus for any element $z\in S$ there exists $t\in S$ such that ${x^n}z={x^{n+1}t}$ and this implies that ${\rm{Im}}f^n={\rm{Im}}f^{n+1}$. Hence again by Proposition \ref{pr3.2} we have the result.\\
The proof of the other iplications is straightforward.

\end{proof}

\begin{proposition} \label{pr3.8}
 Let $A$ be a finitely generated $S$-act over a monoid $S$. Then the following conditions are equivalent:
\item[{\rm (i)}] All factor acts of $A$ are Hopfian and co-Hopfian.
\item[{\rm (ii)}] All factor acts of $A$ are Fitting.
\end{proposition}

\begin{proof}
(i)$\Rightarrow$ (ii) By Proposition \ref{pr3.5}, $A$ is strongly co-Hopfian. Thus if $f\in End(A)$, then by Proposition \ref{pr3.2}, there exists a natural number $n$ such that $ {\rm{Im}}f^n={\rm{Im}}f^{n+1}= {\rm{Im}}f^{2n}$. Define $g:{\rm{Im}}f^n\longrightarrow {\rm{Im}}f^n={\rm{Im}}f^{2n}$ by $g(x)={f^n}(x)$ for any $x\in S$. Since ${A/\ker  f^n}\cong {\rm{Im}}f^n,~ g$ is an epimorphism and by assumption it is a monomorphism. Hence $\ker  f^n\cap \mathcal{K}_{{\rm{Im}} f^{n}}=\Delta_A$, and by Proposition \ref{pr3.2}, $A$ is strongly Hopfian. Now the proof is completed by Theorem \ref{te3.3}.\\
(ii)$\Rightarrow$ (i) Clear by Theorem \ref{te3.3}.
\end{proof}

By the above proposition and Theorem \ref{te3.3}, the next corollaries can be proved.
\begin{corollary} \label{co3.9}
Let $S$ be a monoid. Then the following conditions are equivalent:
\item[{\rm (i)}] All finitely generated $S$-acts are Hopfian and co-Hopfian.
\item[{\rm (ii)}] All finitely generated $S$-acts are Fitting.
\item[{\rm (iii)}] All finitely generated $S$-acts are strongly Hopfian and strongly co-Hopfian.
\end{corollary}

\begin{corollary} \label{co3.10}
Let $S$ be a monoid. Then the following conditions are equivalent:
\item[{\rm (i)}] All cyclic $S$-acts are Hopfian and co-Hopfian.
\item[{\rm (ii)}] All cyclic $S$-acts are Fitting.
\item[{\rm (iii)}] All cyclic $S$-acts are strongly Hopfian and strongly co-Hopfian.
\end{corollary}


~~~~~~~~~~~~~~~~~~~~~~~~~~~~~~~~~~~~





\end{document}